\documentclass[12pt, twoside, leqno]{article}
\usepackage[utf8]{inputenc}
\usepackage[T1]{fontenc}
\usepackage{amsthm}
\usepackage{amssymb}
\usepackage{mathrsfs}
\usepackage{amsmath}
\usepackage{amsfonts}
\usepackage{mathtools}
\usepackage{hyperref}
\usepackage{enumitem}

\usepackage{tikz}

\mathtoolsset{showonlyrefs}

\newtheorem*{Theo}{Theorem}
\newtheorem*{Prop}{Proposition}
\newtheorem{Le}{Lemma}

\DeclareMathOperator{\lcm}{lcm}

\renewcommand{\P}{\mathcal{P}}
\newcommand{\C}{\mathcal{C}}
\newcommand{\A}{\mathcal{A}}

\begin{document}


\baselineskip=17pt


\title{On path partitions of the divisor graph}

\author{Paul Melotti\\
Sorbonne Université\\ 
LPSM\\
4 place Jussieu\\
F-75005 PARIS\\
FRANCE\\
E-mail: paul.melotti@upmc.fr
\and 
Éric Saias\\
Sorbonne Université\\ 
LPSM\\
4 place Jussieu\\
F-75005 PARIS\\
FRANCE\\
E-mail: eric.saias@upmc.fr
}

\date{}

\maketitle

\renewcommand{\thefootnote}{}



\renewcommand{\thefootnote}{\arabic{footnote}}
\setcounter{footnote}{0}

\begin{abstract}
  It is known that the longest simple path in the divisor
  graph that uses integers $\leq N$ is of length $\asymp N/\log N$. We
  study the partitions of $\{1,2,\dots, N\}$ into a \emph{minimal}
  number of paths of the divisor graph, and we show that in such a
  partition, the longest path can have length asymptotically
  $N^{1-o(1)}$.
\end{abstract}

\section{Introduction}
\label{sec:introduction}

The \textit{divisor graph} is the unoriented graph whose vertices are
the positive integers, and edges are the $\{a,b\}$ such that
$a < b$ and $a$ divides $b$. A \textit{path} of length $l$ in the
divisor graph is a finite sequence $n_1,\dots,n_l$ of pairwise
distinct positive integers such that $n_i$ is either a divisor or a
multiple of $n_{i+1}$, for all $i$ such that $1\leq i < l$. Let $F(x)$
be the minimal cardinal of a partition of $\{1,2,\dots,\lfloor x
\rfloor\}$ into paths of the divisor graph.

The asymptotic behaviour of $F(x)$ has been studied in
\cite{ErdosSaias,Saias:div3,Mazet,Chadozeau}. Thanks to the works of
Mazet and Chadozeau, we know that there is a constant
$c\in (\frac16,\frac14)$ such that
\begin{equation}
  \label{eq:1}
  F(x) = c x \left(1 + O\left(\frac{1}{\log \log x \ \log \log \log x}
      \right) \right).
\end{equation}

A partition of $\{1,2,\dots,N\}$ into paths of the divisor graph is
said to be \textit{optimal} if its cardinal is $F(N)$. We are
interested in the length of the paths in an optimal partition.

Let us take the example $N=30$ that was considered in
\cite{Pomerance,Saias:div3}. It is known (see \cite{Saias:div3}) that
$F(30)=5$, so that the following partition is optimal:

\begin{equation}
  \label{eq:2}
  \begin{split}
    & 13, 26, 1, 11, 22, 2, 14, 28, 7, 21, 3, 27, 9, 18, 6, 12, 24, 8,
    16, 4, 20, 10, 30, 15, 5, 25 \\
    & 17 \\
    & 19 \\
    & 23 \\
    & 29
  \end{split}
\end{equation}

Four of these five paths are singletons. In fact, at the end of the
proof of Theorem~2 of \cite{ErdosSaias}, it is proven that the number
of singletons in a (not necessarily optimal) partition is $\asymp N$
for $N$ large enough.

Let us look at the longest paths in an optimal partition of
$\{1,2,\dots,N\}$. Let $L(N)$ be the maximal path length, among all
paths of all \emph{optimal} partitions of $\{1,2,\dots,N\}$ into paths
of the divisor graph. Let also $f(N)$ denote the maximal length of a
path of the divisor graph that uses integers $\leq N$.

It is known that (Theorem~2 of \cite{Saias:dense})
\begin{equation}
  \label{eq:3}
  f(N) \asymp \frac{N}{\log N}.
\end{equation}
Of course $L(N) \leq f(N)$. In the previous example, four of the
five paths are singletons, which implies that the longest path has
maximal length. In other words $L(30) = f(30) = 26$.
More generally, we know that for all $N \geq 1$,
\begin{equation}
  \label{eq:4}
  F(N) \geq N - \lfloor N/2 \rfloor - \lfloor N/3 \rfloor
\end{equation}
(see \cite{Saias:div3}). Inspired by the case
$N=30$, for any $N \in [1,33]$ it is easy to construct a partition of
$\{1,\dots,N\}$ into $N - \lfloor N/2 \rfloor - \lfloor N/3 \rfloor$
paths, all of them but one being singletons. This shows that for $1
\leq N \leq 33$, \eqref{eq:4} is an equality and $L(N) = f(N) = \lfloor
N/2 \rfloor + \lfloor N/3 \rfloor + 1$.


However for larger $N$ the situation becomes more complicated. For $N$
large enough there is no optimal partition with all paths but one
being singletons. This can be deduced from \eqref{eq:3} and the fact
that the constant $c$ in \eqref{eq:1} is less than $1$. Still, it is
natural to wonder if the equality $L(N)=f(N)$ holds for any $N\geq 1$.

We were unable to answer this question, but we looked for lower bounds
on $L(N)$ and proved the following:

\begin{Theo}
  There is a constant $A \geq 0$ such that for all $N\geq 3$,
  \begin{equation}
    \label{eq:5}
    L(N) \geq \frac{N}{(\log N)^A \exp\left[\frac{(\log \log
          N)^2}{\log 2} \right]}.
  \end{equation}
\end{Theo}

To prove this we introduce a new function $H(x)$. For a real number
$x\geq 1$ and two distinct integers $a,b \in [1,x]$, let $L_{a,b}(x)$
be the maximal length of a path having $a$ and $b$ as endpoints and
belonging to an \emph{optimal} partition of $\{1,2,\dots,\lfloor x
\rfloor\}$. If there is no such path, we set $L_{a,b}(x) = 0$.
Then we set
\begin{equation}
  \label{eq:6}
  H(x) = \min L_{r',r}(x)
\end{equation}
where the $\min$ is over all couples $(r',r)$ of \emph{prime} numbers
such that
\begin{equation}
  \label{eq:7}
  \frac{x}{3} < r \leq \frac{x}{2} < r' \leq x.
\end{equation}

The theorem will be an easy consequence of the following.
\begin{Prop}
  There is a constant $N_0$ such that for any $N\geq N_0$, there is a
  set $\P(N)$ of prime numbers in $(3 \sqrt{N \log N}, 4 \sqrt{N \log
    N}]$, of cardinal $|\P(N)| \geq \frac{\sqrt{N}}{19 (\log N)^{3/2}}$, such
  that
  \begin{equation}
    \label{eq:8}
    H(N) \geq \sum_{p \in \P(N)}H\left(\frac{N}{p}\right).
  \end{equation}
\end{Prop}

The technique used here is analogous to that of \cite{Saias:Buchstab}
in the study of the longest path. More precisely, in
\cite{Saias:Buchstab}, $f^*(N)$ denotes the maximal length of a path
that uses integers in $[\sqrt{N},N]$. A quantity $h^*$ is introduced,
which is to $f^*$ what $H$ is to $L$ in our case. The inequality
\eqref{eq:8} is analogous to Buchstab's unequality $(40)$ from
\cite{Saias:Buchstab}. The corresponding lower bounds led to the proof
that $f^*(N) \asymp N/\log N$ (Theorem~2 in \cite{Saias:dense}).

The analogy can be pushed further: in both the proof of \eqref{eq:8}
and of $(40)$ in \cite{Saias:Buchstab}, we borrow a technique used by
Erdős, Freud and Hegyvári who proved the following asymptotic
behaviour:
\begin{equation}
  \label{eq:9}
  \min \max_{1\leq i \leq N-1} \lcm (a_i,a_{i+1}) = \left(\frac14 +
    o(1) \right) \frac{N^2}{\log N},
\end{equation}
where the $\min$ is over all permutations $(a_1,a_2,\dots,a_N)$ of
$\{1,2,\dots,N\}$; see Theorem~1 of \cite{ErdosFreudHegyvari}. In
\cite{ErdosFreudHegyvari} as in \cite{Saias:Buchstab} or in the
present work, the proof goes through the construction of a sequence of
integers by concatenating blocks whose largest prime factor is
constant, and linking blocks together with separating integers. In
\cite{Saias:Buchstab} as in the present work, these blocks take the
form of sub-paths $p \C_{N/p}$, where the $\C_{N/p}$ is a path of
integers $\leq N/p$ whose largest prime factor is $\leq p$.

It is worth mentioning that the article \cite{ErdosFreudHegyvari} of
Erdős, Freud and Hegyvári is the origin of all works related to the
divisor graph.

\section{Notations}
\label{sec:notations}

The letters $p,q,q',r,r'$ will always denote generic prime
numbers. For an integer $m\geq 2$, $P^-(m)$ denotes the smallest prime
factor of $m$.

Let $N\geq 1$. A \emph{path} of integers $\leq N$ of length $l$ is a
$l$-uple $\C = (a_1,a_2,\dots,a_l)$ of pairwise distinct positive
integers $\leq N$, such that for all $i$ with $1\leq i \leq l-1$,
$a_i$ is either a divisor or a multiple of $a_{i+1}$. For convenience,
we take $\C$ up to global flip, \textit{i.e.} we identify
$(a_1,\dots,a_l)$ with $(a_l,\dots,a_1)$. We will denote this path by
$a_1-a_2-\dots-a_l$ (or $a_l-\dots-a_2-a_1$). If $b$ and $c$ are
integers such that $b=a_i$ and $c=a_{i\pm 1}$ for some $i$, we say
that $b$ and $c$ are \emph{neighbours} (in $\C$).

When a partition $\A(N)$ of $\{1,2,\dots,N\}$ is fixed, for any $n \in
\{1,2,\dots,N\}$ we will simply denote by $\C(n)$ the path that
contains $n$ in $\A(N)$.

A partition of $\{1,2,\dots,N\}$ into paths is said to be
\emph{optimal} if it contains $F(N)$ paths (see the
Introduction for the definition of $F$).

Let $\C$ be a path of integers $\leq N$ and $1\leq n \leq N$. Then
$\C$ is said to be $n$-\emph{factorizable} if all the integers of $\C$
are multiple of $n$. Then $\C$ can be written as $\C = n \mathcal{D}$
where $\mathcal{D}$ is a path of integers $\leq N/n$.

For integers $1\leq n \leq N$ and a partition $\mathcal{A}(N)$ of
$\{1,2,\dots,N\}$, we say that $n$ is \emph{factorizing} for
$\mathcal{A}(N)$ if every path of $\mathcal{A}(N)$ that contains a
multiple of $n$ is $n$-factorizable.

\section{Lemmas}
\label{sec:lemmas}

\begin{Le}
  \label{le:1}
  Let $N\geq 1$ and $\A(N)$ be an optimal partition.
  \begin{enumerate}[label=(\roman*)]
  \item
    \label{item:i}
    Let $1\leq n \leq N$ with $n$ factorizing for $\A(N)$. Let $k
    = \lfloor N/n \rfloor$. There are exactly $F(k)$ paths in $\A(N)$
    that contain a multiple of $n$. They are of the form
    $n\mathcal{D}_1, n\mathcal{D}_2, \dots,n\mathcal{D}_{F(k)}$ where
    $\mathcal{D}_1, \mathcal{D}_2, \dots,\mathcal{D}_{F(k)}$ is an
    optimal partition of $\{1,2,\dots,k\}$.
  \item
    \label{item:ii}
    Let $z>1$ be a real number. Let $M_z(N)$ be the set of
    integers $m\leq N$ that are \emph{not} factorizing for $\A(N)$ and such that
    \begin{equation}
      \label{eq:10}
      m> \frac{N}{z} \ \text{ and } \ P^-(m)>z.
    \end{equation}
    Then
    \begin{equation}
      \label{eq:11}
      |M_z(N)|<\frac{2N}{z}.
    \end{equation}
  \end{enumerate}
\end{Le}

\begin{proof}
  \begin{enumerate}[label=(\roman*)]
  \item
    The set of paths that contain a multiple of $n$ is of the form
    $\{n\mathcal{D}_1,n\mathcal{D}_2,\dots,n\mathcal{D}_g\}$ where
    $\mathcal{D}_1,\mathcal{D}_2,\dots,\mathcal{D}_g$ is a partition
    of the integers $\leq k = \lfloor N/n \rfloor$. Since $\A(N)$ is
    optimal, $\mathcal{D}_1,\mathcal{D}_2,\dots,\mathcal{D}_g$ is
    optimal, hence $g=F(k)$.
  \item Let $m\in M_z(N)$. There is a path $\C_m$ in $\A(N)$ that
    contains multiples and non-multiples of $m$. Hence there is an
    integer $c(m)$ in $\C_m$ that is not a multiple of $m$, and is
    neighbour to an integer $b(m)$ which is a multiple of $m$. Then
    $c(m)$ has to be a divisor of $b(m)$. More precisely, if
    $b(m) = a m$, then $c(m)$ can be written as $c(m)=\tilde{a} \tilde{m}$
    with $\tilde{a}$ a divisor of $a$ and $\tilde{m}$ a strict divisor of
    $m$. Since $P^-(m)>z$, $c(m) < N/z$.

    Moreover, if $m,m'$ are two distinct elements of $M_z(N)$, then
    \begin{equation}
      \label{eq:12}
      \lcm(m,m') \geq \min(m P^-(m'), m' P^-(m)) > \frac{N}{z} z = N.
    \end{equation}
    As a result the map
    \begin{equation}
      \label{eq:13}
      \begin{split}
        b: M_z(N) & \to \{1,2,\dots,N\} \\
        m & \mapsto b(m)
      \end{split}
    \end{equation}
    is an injection.

    Moreover, any integer $c < N/z$ has at most two neighbours in
    $\C(c)$. Consequently the map
    \begin{equation}
      \label{eq:14}
      \begin{split}
        c: M_z(N) & \to \{1\leq n < N/z \} \\
        m & \mapsto c(m)
      \end{split}
    \end{equation}
    is at-most-two-to-one. Thus
    \begin{equation}
      \label{eq:14b}
      |M_z(N)| < \frac{2N}{z}.
    \end{equation}
  \end{enumerate}
\end{proof}

\begin{Le}
  \label{le:2}
  There exists a constant $N_1$ such that for any $N \geq N_1$, there
  is a set $\widetilde{\mathcal{P}}(N)$ of prime numbers in $(3\sqrt{N
    \log N}, 4\sqrt{N \log N}]$ of cardinal
  \begin{equation}
    \label{eq:sizePt}
    |\widetilde{\mathcal{P}}(N)| \geq \sqrt{\frac{N}{\log N}},
  \end{equation}
  such that for any prime numbers $r,r'$ with
  \begin{equation}
    \label{eq:15}
    \frac{N}{3} < r \leq \frac{N}{2} < r' \leq N,
  \end{equation}
  there exists an optimal partition $\A(N)$ of $\{1,2,\dots,N\}$ that
  contains the paths $r'$ and $2r-r$ and for which all the integers in
  $\widetilde{\mathcal{P}}(N)$ are factorizing.
\end{Le}

\begin{proof}
  Let $N_1$ be such that for any $N\geq N_1$,
  \begin{align}
    \label{eq:16}
    & \pi\left(4\sqrt{N\log N}\right) - \pi\left(3\sqrt{N\log N}\right) - \frac23
      \sqrt{\frac{N}{\log N}} \geq \sqrt{\frac{N}{\log N}}, \\
    \label{eq:16b}
    & \pi\left(\frac{N}{2}\right) - \pi\left(\frac{N}{3}\right) \geq 8.
  \end{align}
  The existence of such a $N_1$ comes from the prime number theorem
  (more precisely the left-hand-side of \eqref{eq:16} is equivalent to
  $\frac43 \sqrt{\frac{N}{\log N}}$). We also take $N_1$ large enough
  so that
  \begin{equation}
    \label{eq:inter}
    \left(3\sqrt{N\log N}, 4 \sqrt{N \log N}\right] \cap
    \left(\frac{N}{3},\frac{N}{2}\right] = \emptyset.
  \end{equation}

  Let $N\geq N_1$. We start by fixing an optimal partition
  $\A'(N)$. We apply Lemma~\ref{le:1} \ref{item:ii} to $\A'(N)$ with
  $z=3\sqrt{N \log N}$. All the prime numbers $p$ in
  $(3\sqrt{N\log N}, 4\sqrt{N \log N}]$ that are not factorizing are
  in $M_z(N)$, since they satisfy
  $p > 3\sqrt{N \log N} \geq \frac{N}{z}$ and $P^-(p) = p > z$, so
  there are at most $\frac23 \sqrt{\frac{N}{\log N}}$ of them. By
  removing these and using \eqref{eq:16}, we get a set
  $\widetilde{\mathcal{P}}(N)$ of prime numbers in
  $(3\sqrt{N\log N}, 4\sqrt{N \log N}]$ that are factorizing in
  $\A'(N)$, with cardinality
  \begin{equation}
    \label{eq:17}
    |\widetilde{\mathcal{P}}(N)| \geq \sqrt{\frac{N}{\log N}}.
  \end{equation}

  We now change notations slightly and fix two prime numbers
  $r_0, r'_0$ such that
  \begin{equation}
    \label{eq:18a}
    \frac{N}{3} < r_0 \leq \frac{N}{2} < r'_0 \leq N.
  \end{equation}
  Our goal is to go from $\A'(N)$ to a new optimal partition $\A(N)$
  that contains the paths $r'_0$ and $2r_0-r_0$ \emph{while
    maintaining the fact that the elements of
  }$\widetilde{\mathcal{P}}(N)$ \emph{are factorizing}.

  Let us denote the set of prime numbers
  \begin{equation}
    \label{eq:24}
    \mathcal{R} = \left\{ \frac{N}{3} < r \leq \frac{N}{2} \right\},
  \end{equation}
  and $\mathcal{R}^*(\A'(N))$ the subset of
  $r \in \mathcal{R}$ such that $r$ does not have $1$ as a neighbour
  in $\C(r)$ and $2r$ does not have $1$ nor $2$ has a neighbour in
  $\C(2r)$. Then for any $r\in \mathcal{R}^*(\A'(N))$, since the only possible
  neighbour of $r$ is $2r$ and reciprocally, by optimality the path
  $\C(r)$ is equal to $r-2r$. Moreover, since $1$ and $2$ have at most two
  neighbours,
  \begin{equation}
    \label{eq:18}
    |\mathcal{R} \setminus \mathcal{R}^*(\A'(N))| \leq 4.
  \end{equation}

  Now we make it so that $r'_0$ is a path. If it is not the case,
  since the only possible neighbour of $r'_0$ is $1$, $\C(r'_0)$ is of
  the form $\mathcal{D}-r'_0$ with $\mathcal{D}$ a path ending in
  $1$. We split this path into $\mathcal{D}$ on one side and $r'_0$ on
  the other side. By \eqref{eq:18} and \eqref{eq:16b}, there is at
  least one element $r^*\in \mathcal{R}^*(\A'(N))$. We stick
  $\mathcal{D}$ to $\C(r^*)$, thus forming the path
  $\mathcal{D}-\C(r^*)$. This is possible because $\mathcal{D}$ ends in
  $1$. Let $\A''(N)$ be this new partition. The total number of paths
  has not changed so $\A'(N)$ is still optimal, furthermore it
  contains the path $r'_0$, and the elements of
  $\widetilde{\mathcal{P}}(N)$ are still factorizing because the
  integers in the paths
  that changed were not multiples of any
  $p\in \widetilde{\mathcal{P}}(N)$.

  The subset $\mathcal{R}^*(\A''(N))$ might differ from
  $\mathcal{R}^*(\A'(N))$ by one element, but it still satisfies
  \eqref{eq:18} and its elements $r$ still satisfy that $\C(r)$ is
  equal to $r-2r$. If $r_0 \in \mathcal{R}^*(\A''(N))$, we can set
  $\A(N) = \A''(N)$ and the proof is over. We now suppose that
  $r_0 \notin \mathcal{R}^*(\A''(N))$.

  By \eqref{eq:18} and \eqref{eq:16b}, there are at least four
  elements $r_1,r_2,r_3,r_4$ in $\mathcal{R}^*(\A''(N))$. We cut the
  path $\C(1)$ into one, two or three paths, one of them being the
  singleton $1$ (we will see later that we get in fact three
  paths). Such a move will be called an \emph{extraction} of the
  integer $1$. We similarly \emph{extract} the integer $2$. We now use
  these integers $1$ and $2$ to stick together the paths $r_i-2r_i$ by
  forming
  \begin{equation*}
    r_1-2r_1-1-2r_2-r_2 \ \text{ and } \ r_3-2r_3-2-2r_4-r_4.
  \end{equation*}
  We thus get a new partition $\A(N)$. Its number of paths is less or
  equal to that of $\A''(N)$, so it is still optimal (this shows in
  particular that $1$ and $2$ were not endpoints of their paths). It
  also satisfies $r_0 \in \mathcal{R}^*(\A(N))$ since $1$ and $2$ are
  not linked to $r_0$ nor $2r_0$, so that it contains the path
  $r_0-2r_0$, as well as $r'_0$, and the elements of
  $\widetilde{\mathcal{P}}(N)$ are still factorizing.
\end{proof}

\section{Proof of the Proposition}
\label{sec:proof-proposition}


Let $N_1$ be the constant of Lemma~\ref{le:2}. We fix a $N_0$ such
that
\begin{equation}
  \label{eq:19}
  N_0 \geq N_1^4
\end{equation}
and such that for all $N\geq N_0$,
\begin{equation}
  \label{eq:20}
  \begin{split}
    \frac12 \sqrt{\frac{N}{\log N}}
   & \geq \pi\left(\frac14
    \sqrt{\frac{N}{\log N}} \right) - \pi\left(\frac16
    \sqrt{\frac{N}{\log N}}\right) \\
   & \ \ \geq \pi\left(\frac18
    \sqrt{\frac{N}{\log N}} \right) - \pi\left(\frac19
    \sqrt{\frac{N}{\log N}}\right) \\
   & \ \ \ \ \geq \left\lfloor \frac{\sqrt{N}}{37 (\log N)^{3/2}} \right\rfloor
    \geq \frac{\sqrt{N}}{38 (\log N)^{3/2}} + \frac12
    \geq 5
  \end{split}
\end{equation}
and
\begin{equation}
  \label{eq:21}
  4 \sqrt{\log N} \leq N^{1/4}.
\end{equation}
The existence of such a $N_0$ is again an easy consequence of the
prime number theorem. Also note that since $N_0 \geq N_1$,
\eqref{eq:inter} still holds.

Let $N\geq N_0$.
We chose a set $\widetilde{\P}(N)$ according to Lemma~\ref{le:2}.
Let us denote
\begin{equation}
  \label{eq:I}
  I = \left\lfloor \frac{1}{37} \frac{\sqrt{N}}{(\log N)^{3/2}} \right\rfloor.
\end{equation}
By \eqref{eq:sizePt} and \eqref{eq:20} we can chose $2I$ elements in
$\widetilde{\P}(N)$, which we denote as
\begin{equation}
  \label{eq:27}
  p_1,p_2,\dots,p_{2I}.
\end{equation}
We set $\P(N) = \{ p_1,\dots,p_{2I-1}\}$.  By \eqref{eq:20} again,
$|\mathcal{P}(N)| \geq \frac{\sqrt{N}}{19 (\log N)^{3/2}}$. It remains
to prove that this set $\P(N)$ satisfies \eqref{eq:8}.

Let $r,r'$ be two prime numbers such that
\begin{equation}
  \label{eq:22}
  \frac{N}{3} < r \leq \frac{N}{2} < r' \leq N.
\end{equation}
By the property of $\widetilde{\P}(N)$ in Lemma~\ref{le:2}, there
exists an optimal partition $\A'(N)$, that contains the paths $r'$ and
$2r-r$, for which the elements of $\widetilde{\mathcal{P}}(N)$ (and in
particular the elements of $\P(N)$) are
factorizing.

We denote two sets of prime numbers
\begin{equation}
  \label{eq:23}
  \begin{split}
    \mathcal{Q}(N) = \left\{ \frac19 \sqrt{\frac{N}{\log N}} < q \leq
    \frac18 \sqrt{\frac{N}{\log N}} \right\}, \\
    \mathcal{Q}'(N) = \left\{ \frac16 \sqrt{\frac{N}{\log N}} < q' \leq
    \frac14 \sqrt{\frac{N}{\log N}} \right\}. 
  \end{split}
\end{equation}
For all $(p,q,q') \in \widetilde{\mathcal{P}}(N) \times \mathcal{Q}(N)
\times \mathcal{Q}'(N)$ we have
\begin{align}
  \label{eq:25}
  \frac{N}{3} < pq \leq \frac{N}{2}, \\
  \label{eq:26}
  \frac{N}{2} < pq' \leq N.
\end{align}

We focus on the factorizing prime number $p_{2I}$. For any
$q \in \mathcal{Q}(N)$, because of \eqref{eq:25} the only possible
neighbours of $p_{2I}q$ are $p_{2I}$ and $2p_{2I}q$. Similarly, the
only possible neighbours of $2p_{2I}q$ are $p_{2I}, 2p_{2I}$ or
$p_{2I}q$. But $p_{2I}$ and $2p_{2I}$ can be linked to at most $4$
elements of type $p_{2I}q$ or $2p_{2I}q$. By \eqref{eq:20} we know
that $|\mathcal{Q}(N)| \geq 5$, so there exists a
$q_{2I}\in\mathcal{Q}(N)$ for which neither $p_{2I}q_{2I}$ nor
$2p_{2I}q_{2I}$ is a neighbour of $p_{2I}$ or $2p_{2I}$. As a result,
the only possible neighbour for $p_{2I}q_{2I}$ is $2p_{2I}q_{2I}$, and
reciprocally. By optimality, $\A'(N)$ contains the path
$p_{2I}q_{2I} - 2p_{2I}q_{2I}$.

Using \eqref{eq:20} we can chose
\begin{itemize}  
\item
  $I$ elements of $\mathcal{Q}'(N)$ which we write as
  \begin{equation}
    \label{eq:28}
    q_1,q_3,\dots,q_{2I-1};
  \end{equation}
\item
  $I-1$ elements of $\mathcal{Q}(N) \setminus \{q_{2I}\}$ which we write
  as
  \begin{equation}
    \label{eq:28b}
    q_2,q_4,\dots,q_{2I-2}.
  \end{equation}
\end{itemize}

Let $i$ be such that $1\leq i \leq 2I-1$. Then the prime number $p_i$
is factorizing for $\A'(N)$ so by Lemma~\ref{le:1}~\ref{item:i} the
paths of $\A'(N)$ that contain multiples of $p_i$ are of the form
\begin{equation}
  \label{eq:29}
  p_i \C_{i,1}, p_i \C_{i,2}, \dots, p_i\C_{i,F(N/p_i)}
\end{equation}
where $ \C_{i,1},\C_{i,2}, \dots, \C_{i,F(N/p_i)}$ is an optimal
partition of $\{1,2,\dots,\lfloor N/p_i \rfloor\}$. By our choice of
indices \eqref{eq:28},\eqref{eq:28b}, one of the elements
$q_i,q_{i+1}$ is in $\mathcal{Q}'(N)$, we rename it $\widetilde{q_i}$,
and the other is in $\mathcal{Q}(N)$, we rename it
$\widetilde{q_{i+1}}$. Using \eqref{eq:25},\eqref{eq:26} we get
\begin{equation}
  \label{eq:30}
  \frac{N}{3p_i} < \widetilde{q_{i+1}} \leq \frac{N}{2p_i} <
  \widetilde{q_{i}} \leq \frac{N}{p_i}.
\end{equation}
Using \eqref{eq:21} and \eqref{eq:19}, we have
$N/p_i \geq N^{1/4} \geq N_0^{1/4} \geq N_1$. Hence we can apply
Lemma~\ref{le:2} with $N/p_i$ instead of $N$. We deduce that there
exists an optimal partition of $\{1,2,\dots,\lfloor N/p_i\rfloor\}$
that contains the paths $\widetilde{q_i}$ and
$\widetilde{q_{i+1}} - 2\widetilde{q_{i+1}}$.  By extracting $1$ in
that partition, we can stick these two paths together into
$\widetilde{q_i}-1-2\widetilde{q_{i+1}} - \widetilde{q_{i+1}}$ while
keeping an optimal partition. To sum up, we know now that there is an
optimal partition of the integers $\leq N/p_i$ containing a path that
has $q_i$ and $q_{i+1}$ as endpoints.

Let
$\mathcal{D}_{i,1}, \mathcal{D}_{i,2},\dots, \mathcal{D}_{i,F(N/p_i)}$
be an optimal partition of the integers $\leq N/p_i$, with
$\mathcal{D}_{i,1}$ having $q_i,q_{i+1}$ as endpoints and of maximal
length $L_{q_i,q_{i+1}}(N/p_i)$.  We can transform $\A'(N)$ by
replacing the paths $\left(p_i\C_{i,j}\right)_{1\leq j \leq F(N/p_i)}$
by $\left(p_i\mathcal{D}_{i,j}\right)_{1\leq j \leq F(N/p_i)}$. In
this way we get a new optimal partition $\A''(N)$ that contains all
the paths $p_i\mathcal{D}_{i,1}$ for $1\leq i \leq 2I-1$, as well as
$r'$, $2r-r$, and $p_{2I}q_{2I} - 2p_{2I}q_{2I}$.

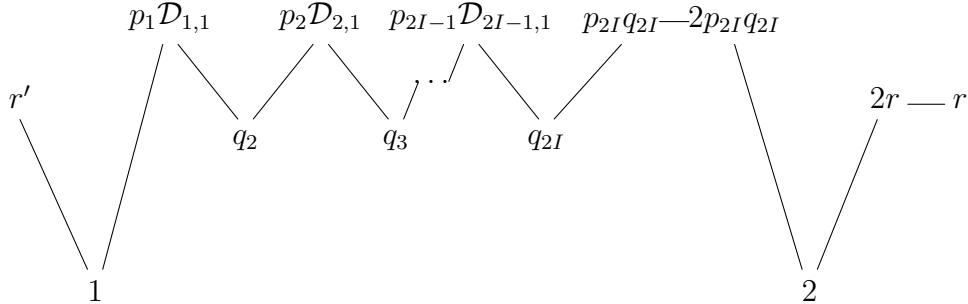
\begin{figure}[h]
  \centering
  \begin{tikzpicture}
    \draw (-1,1) node [above] {$r'$};
    \draw (0,-1) node [below] {$1$};
    \draw (1,2) node [above] {$p_1\mathcal{D}_{1,1}$};
    \draw (2,1) node [below] {$q_2$};
    \draw (3,2) node [above] {$p_2\mathcal{D}_{2,1}$};
    \draw (4,1) node [below] {$q_3$};
    \draw (4.5,1.5) node {$\dots$};
    \draw (5,2) node [above] {$p_{2I-1} \mathcal{D}_{2I-1,1}$};
    \draw (6,1) node [below] {$q_{2I}$};
    \draw (7,2) node [above] {$p_{2I}q_{2I}$};
    \draw (8.5,2) node [above] {$2p_{2I}q_{2I}$};
    \draw (9.5,-1) node [below] {$2$};
    \draw (10.5,1) node [above] {$2r$};
    \draw (11.5,1) node [above] {$r$};
    \draw (-1,1) -- (-0.1,-1);
    \draw (0.1,-1) -- (0.9,2);
    \draw (1.1,2) -- (1.9,1);
    \draw (2.1,1) -- (2.9,2);
    \draw (3.1,2) -- (3.9,1);
    \draw (4.1,1) -- (4.3,1.5);
    \draw (4.7,1.5) -- (4.9,2);
    \draw (5.1,2) -- (5.9,1);
    \draw (6.1,1) -- (7,2);
    \draw (7.5,2.3) -- (7.9,2.3);
    \draw (8.5,2) -- (9.4,-1);
    \draw (9.6,-1) -- (10.4,1);
    \draw (10.8,1.2) -- (11.3,1.2);
  \end{tikzpicture}
  \caption{A long path with endpoints $r',r$.}
  \label{fig}
\end{figure}

By extracting the integers $1$, $2$ and the $q_i$ for $2\leq i \leq
2I$, we construct the path of Figure~\ref{fig} while keeping an optimal
partition of $\{1,2,\dots,N\}$.

Its length is larger than
\begin{equation}
  \label{eq:31}
    \sum_{i=1}^{2I-1}L_{q_i,q_{i+1}}(N/p_i)
    \geq \sum_{p\in \mathcal{P}(N)} H(N/p_i).
\end{equation}

This being true for any $r,r'$ satisfying \eqref{eq:22}, we get

\begin{equation}
  \label{eq:32}
  H(N) \geq \sum_{p\in \mathcal{P}(N)} H(N/p_i).
\end{equation}
\qed

\section{Proof of the Theorem}
\label{sec:proof-theorem}

Let us fix a constant $N_2 = 2^{2^{k_0}} \geq N_0$, where $N_0$ is the
constant from the Proposition.
We chose a constant $B$ such that for all $N \leq 2^{2^{k_0+2}}$,
\begin{equation}
  \label{eq:33}
  N \leq 4 (\log N)^B \exp\left[\frac{(\log \log N)^2}{\log 2}\right]
\end{equation}
and
\begin{equation}
  \label{eq:34a}
  B \geq 8.
\end{equation}

We show by induction on $k \geq k_0 + 2$ that for all $N$ such that
\begin{equation}
  \label{eq:35}
  2^{2^{k_0}} < N \leq 2^{2^k},
\end{equation}
 we have
\begin{equation}
  \label{eq:34}
  H(N) \geq \frac{N}{(\log N)^B \exp\left[\frac{(\log \log N)^2}{\log
        2} \right]}
\end{equation}

\subsubsection*{Base case}

Let $N$ be such that $2^{2^{k_0}} < N \leq 2^{2^{k_0 + 2}}$, then we
have $N > N_2 \geq N_0 \geq N_1^4$ (see \eqref{eq:19}) with $N_1$ the
constant of Lemma~\ref{le:2}. Let $r,r'$ be two prime numbers such
\begin{equation}
  \label{eq:36}
  \frac{N}{3} < r \leq \frac{N}{2} < r' \leq N.
\end{equation}

Lemma~\ref{le:2} implies that there is an optimal partition $\A(N)$ of
$\{1,2,\dots,N\}$ which contains the paths $r'$ and $2r-r$. By
extracting $1$, we can stick them into $r'-1-2r-r$ while keeping an
optimal partition. This implies that $H(N) \geq 4$, and \eqref{eq:33}
yields the base case.

\subsubsection*{Induction step}

Let $k \geq k_0+2$. We suppose that \eqref{eq:34} holds for all $N \in
\left(2^{2^{k_0}},2^{2^k}\right]$.

Let $N$ be such that $2^{2^{k}} < N \leq 2^{2^{k+1}}$.
Since $k\geq k_0+2$, we also have $N^{1/4} > 2^{2^{k_0}}$.

Let $p \in \left(3\sqrt{N\log N}, 4\sqrt{N \log N} \right]$. By
\eqref{eq:21}, we have
\begin{equation}
  \label{eq:37}
  2^{2^{k_0}} < N^{1/4} \leq \frac{N}{p} \leq \sqrt{N} \leq 2^{2^k}.
\end{equation}
By using the induction hypothesis on $N/p$, we get

\begin{align}
  H\left(\frac{N}{p}\right) & \geq \frac{N}{p(\log(N/p))^B \exp\left[\frac{(\log \log (N/p))^2}{\log
                              2} \right]} \\
  & \geq \frac{N}{p(\log\sqrt{N})^B \exp\left[\frac{(\log \log \sqrt{N})^2}{\log
    2} \right]} \\
  & = \frac{2^{B-1} (\log N)^2 N}{p(\log N)^B \exp\left[\frac{(\log
    \log N )^2}{\log 2} \right]}.
\end{align}

Hence by using the Proposition and \eqref{eq:34a},
\begin{align}
  H(N) & \geq \sum_{p \in \mathcal{P}(N)} H\left(\frac{N}{p}\right) \\
  & \geq \frac{|\mathcal{P}(N)|}{\max \mathcal{P}(N)} \  \frac{2^{B-1}
    (\log N)^2 N}{(\log N)^B \exp\left[\frac{(\log \log N )^2}{\log 2}
    \right]}\\
  & \geq \frac{2^{B-1}}{76} \frac{N}{ (\log N)^B \exp\left[\frac{(\log
    \log N )^2}{\log 2} \right]} \\
  & \geq \frac{N}{ (\log N)^B \exp\left[\frac{(\log
    \log N )^2}{\log 2} \right]}.
\end{align}

This concludes the induction step.

\medskip

Finally, since $L(N) \geq 1$ for all $N\geq 1$, we get the Theorem by
chosing $A=\max(B,A_0)$ where $A_0$ is a constant such that for all
$3\leq N < N_0$,
\begin{equation}
  \label{eq:38}
  N \leq (\log N)^{A_0} \exp\left[\frac{(\log \log N)^2}{\log 2}\right].
\end{equation}
\qed


\bibliographystyle{acm}
\bibliography{chaines}

\end{document}